\documentclass[a4paper]{article}
\usepackage[margin=1in]{geometry} 
\usepackage[utf8]{inputenc}
\usepackage{amsmath}
\usepackage{array}
\usepackage{amsthm}
\usepackage{amssymb}
\usepackage{graphicx}
\usepackage{algorithmic}
\usepackage{algorithm}
\usepackage[table]{xcolor}
\usepackage{subcaption}
\usepackage{multicol}
\usepackage{multirow}
\usepackage[utf8]{inputenc}
\usepackage{graphicx, color}
\usepackage{tabularx}
\usepackage[pagebackref=true,breaklinks=true,letterpaper=true,colorlinks,bookmarks=false]{hyperref}
\usepackage{xifthen}
\usepackage{mathtools}
\usepackage{url}
\usepackage{bm}
\usepackage{tikz} 
\NewDocumentCommand{\vect}{ O{} O{} m }{\bm{#3}\ifthenelse{\isempty{#1}}{}{^{(#1)}}\ifthenelse{\isempty{#2}}{}{_{#2}}}

\NewDocumentCommand{\mat}{ O{} O{} m }{\bm{#3}\ifthenelse{\isempty{#1}}{}{^{(#1)}}\ifthenelse{\isempty{#2}}{}{_{#2}}}

\usepackage{mathrsfs} 
\newcommand{\comment}[1]{}

\title{Matrix Diagonalization as a Board Game: Teaching an Eigensolver the Fastest Path to Solution}
\author{Phil Romero$^1$\thanks{Corresponding author: \texttt{prr@lanl.gov}} \and Manish Bhattarai$^2$ \and Christian F. A. Negre$^2$ \and Anders M. N. Niklasson$^2$ \and Adetokunbo Adedoyin$^3$}

\date{
    $^1$High Performance Computing division, Los Alamos National Laboratory, NM, USA \\ %
    $^2$Theoretical division, Los Alamos National Laboratory, NM, USA \\ %
    $^3$ Computer, Computational, and Statistical Sciences Division, Los Alamos National Laboratory, Los Alamos, New Mexico 87545
}

\begin{document}
    \maketitle
    
    \begin{abstract}
        Matrix diagonalization is at the cornerstone of numerous fields of scientific computing. Diagonalizing a matrix to solve an eigenvalue problem requires a sequential path of iterations that eventually reaches a sufficiently converged and accurate solution for all the eigenvalues and eigenvectors. This typically translates into a high computational cost. Here we demonstrate how reinforcement learning, using the AlphaZero framework, can accelerate Jacobi matrix diagonalizations by viewing the selection of the fastest path to solution as a board game. To demonstrate the viability of our approach we apply the Jacobi diagonalization algorithm to symmetric Hamiltonian matrices that appear in quantum chemistry calculations. We find that a significant acceleration can often be achieved.  Our findings highlight the opportunity to use machine learning as a promising tool to improve the performance of numerical linear algebra. 
        
    \end{abstract}
    
    
    \section{Introduction}
    \label{introduction}
    Matrix diagonalization is a fundamental procedure in numerical linear algebra and appears in numerous fields of scientific computing. 
    A symmetric matrix $A$, satisfying $A = A^t$ can be factorized as $A = U^tDU$, where $D$ and $U$ are, respectively, matrices containing eigenvalues and orthogonal eigenvectors of $A$ \cite{Watkins2004-qh}. 
    The computational complexity of a matrix diagonalization typically scales cubically, $\mathcal{O}(N^3)$, with the number of eigenpairs, $N$, and often with a fairly large pre-factor. Matrix diagonalizations therefore often become a computational bottleneck. 
    
    Diagonalization is at the core of computational chemistry and it appears as the computational bottleneck determining the properties of chemical systems. In a broad range of quantum chemistry methods based on, for example, Hartree-Fock and density functional theory or semi-empirical methods, a symmetric effective single-particle Hamiltonian, $H$, needs to be diagonalized to construct the density matrix, $\rho$, which is given by the Fermi function, $f$, of $H$, i.e.\ $\rho = f(H)$ \cite{Szabo2012-bc}.
    The decomposition $H=U^t D U$ gives an equivalent system where $D$ is the diagonal matrix containing the energy eigenvalues for each molecular orbital, and $U$ is a unitary matrix containing the eigenvectors corresponding to the molecular orbitals for the system. 
    Hamiltonian matrix functions, such as the density matrix $\rho = f(H)$, can then be constructed in the diagonal molecular-orbital eigenbasis such that $\rho = f(H) = U^t f(D) U$. The Fermi function, $f(H)$, projects the eigenvalues of the occupied states (at lower energies) to 1 and all the unoccupied eigenvalues (at higher energies) to 0.  A property or an observable of a system, ``$a$", such as charges, partial energies and forces, or the dipole moments, can then be determined by computing the trace (or the partial trace) over the corresponding operator, $A$, projected into the occupied subspace by the density matrix, where $a  = \langle A \rangle = tr(\rho A)$. Compared to the prior diagonalization these trace operations are straightforward and can be performed with little computational overhead. 
    
    One of the most extensively used methods for the diagonalization of a matrix $A$ is based on the QR decomposition. The QR method for matrix diagonalization is an iterative solver that uses a sequence of Gram-Schmidt orthogonalization. In the first iteration we perform a QR factorization of $A \equiv A_1 = Q_1 R_1$. We then successively compute $A_{k+1} = R_k Q_k = Q^t_k A_{k} Q_k$, which eventually converges to a sufficiently diagonal matrix $D = Q^t_k Q^t_{k-1} .. Q_0^t A Q_0 Q_{1} .. Q_k = U^t A U$. The diagonal matrix $D$ contains the eigenvalues and the accumulated orthogonal matrix $U$ contains the corresponding eigenvectors.
    This method is surprisingly efficient, but has several limitation. For example, it is not possible to take advantage of data matrix sparsity or data locality to reduce the computational complexity or to achieve efficient parallelism. 
    
    For the special cases of real symmetric matrices, the Jacobi matrix diagonalization method provides an alternative \cite{Rutishauser1966-pe}. The Jacobi method is based on a sequence of Givens rotation matrices, $G^{(i,j)}_k$, that each sets a pair chosen off-diagonal elements, $A_{ij}$ and $A_{ji}$ in $A$ to zero. In the first iteration $A_1 = (G_0^{(i,j)})^t A G_0^{(i,j)}$ for some chosen off-diagonal matrix pair $(i,j)$. The procedure is repeated with $A_{k+1} = (G_{k}^{(i,j)})^t A_k G_{k}^{(i,j)}$.
    Even if other elements besides $A_{ij}$ and $A_{ji}$ are modified in each rotation, the magnitude of the off-diagonal elements are systematically reduced. By iteratively sweeping over non-zero off-diagonal elements using the Givens rotations, the sequence eventually converges to a sufficiently converged diagonal matrix $D = G^t_{k+1} G^t_{k} .. G_0 A_k G_{0} .. G_{k+1} = U^t A U$. The diagonal matrix $D$ contains the eigenvalues and the accumulated orthogonal matrix $U$ contains the corresponding eigenvectors.
    The convergence of the Jacobi diagonalization method depends on the sequence we chose for the Givens rotations. We may, for example, always chose the pairs $(i,j)$ corresponding to the off-diagonal matrix elements with the largest magnitude or we may simply sweep over rows and  columns. Any path can be chosen, but certain paths will converge faster.  
    
    To speed up the matrix diagonalization various methods have been developed that can take advantage of the particular structures of the matrix, such as the matrix symmetry and sparsity~\cite{muro2019acceleration}. Special methods have also been developed when only a few eigenpairs are needed or to take advantage of the available computer architectures, e.g.\ to improve parallelism \cite{vecharynski2015projected}. In this paper we will introduce another aspect that can improve matrix diagonalization, which is the use of {\em machine learning} (ML) to accelerate algorithmic convergence. Our goal is to show how modern machine learning can be used to ``teach" an eigensolver the fastest path to solution to speed up the calculations. The ability to use AI to improve numerical linear algebra is a new promising field of research with some promising recent results \cite{Fawzi2022-qi}.
    
    \section{Methods}
    
    Typically, ML techniques have employed eigensolvers to reduce the dimensionality of the problem to the principal components, hence, eliminating redundant variables~\cite{abdi2010principal}. Here, we do the opposite, which is to use ML and specifically reinforcement learning (RL) to accelerate the computation of the eigenspace of a symmetric matrix. 
    To the extent of our knowledge, no algorithm has yet explored using ML/RL to speedup such calculations. Some recent attempts are available, however, with iterative methods using ML assisted optimizers where no time-to-solution is reported \cite{Han2020-qp}.
    Generalized convolutional deep neural networks have also been proposed to solve the electronic structure problem in quantum chemistry, but without doing an explicit matrix diagonalizaion \cite{Finkelstein2021-sc,Finkelstein2021-xi,Finkelstein2022-xu}.
    
    The recent breakthrough of AlphaZero~\cite{alphago-2016,silver2017mastering}, an RL-based framework, has revolutionized the AI industry by providing solutions for intractable problems in various applications, including protein folding~\cite{alphafold}, and games like Go~\cite{alphago-2016}, Chess, Shogu, and Star Craft \cite{schrittwieser2020mastering}. The use of Monte Carlo tree search (MCTS) and one-step lookahead DNN has allowed AlphaZero to provide heuristic-free exploration. Based on the principle of AlphaZero, a recent work from Google Deep Mind called AlphaTensor demonstrated a proof of concept showing that the existing fastest matrix-matrix multiplication being performed by the Strassen algorithm was significantly accelerated using AI methods \cite{Fawzi2022-qi}. In the aforementioned work, the authors extended the capabilities of the AlphaZero framework to estimate the heuristics for the fastest matrix multiplication.
    
    Jacobi diagonalization method gives an ideal framework to explore ML techniques. Givens rotations are typically constructed using the indices of the maximum absolute off-diagonal elements as a pivot \cite{press1992numerical}. However, exploring a reinforcement learning algorithm could enable the AI to discover unconventional strategies for selecting optimal pivot points and facilitating diagonalization in fewer steps, without prior knowledge of the data generation process or diagonalization strategy.
    The key observation behind this article is that we can view the Jacobi algorithm as a board game. In each move we select a pair of off-diagonal matrix elements, $(i,j)$ and $(j,i)$, on an $N \times N$ board that are removed, i.e.\ that are set to zero by a Givens rotation. Each move leads to an increase in the magnitude of some of the other matrix elements, which can be seen as the counter move by an opponent. The player who can find the fewest number of moves that reaches a sufficiently converged and accurate solution has won. This is a well-defined problem for reinforcement learning, which recently has demonstrated a spectacular breakthrough for board games such as Go and Chess using the framework  \cite{alphago-2016,schrittwieser2020mastering}. 
    Our target problem is the diagonalization of a sequence of symmetric Hamiltonian matrices that appear in quantum-mechanical molecular dynamics (QMD) simulations \cite{Marx2009-zg}. Often tens-of-thousands of fairly similar Hamiltonians need to be diagonalized during a molecular dynamics simulations. This problem should be well-suited for machine learning.
    
    In this work, we present Alpha-FastEigen (FastEigen for short notation), an AlphaZero-based matrix diagonalization framework that achieves accelerated eigen decomposition. FastEigen uniquely learns from the matrix state spaces by combining Monte-Carlo Tree Search (MCTS) and policy-value neural network iteratively to optimize policy and estimate the path with the least number of steps leading to the solution eigenspace. By using a self-play strategy over many episodes, the agent learns a policy that achieves an even faster diagonalization rate. 
    
    \subsection{Data set generation}
    The data set used in this work consists of symmetric Hamiltonian matrices generated from a molecular dynamics trajectory using the LATTE code  \cite{Bock2008-bp}. This is a QMD simulation code where the Hamiltonian matrix is computed according to the Density Functional Tight-Binding (DFTB) theory \cite{Elstner1998-kp}.   
    We created a python library that allows for a quick construction of the Hamiltonian $H$
    from a system coordinate. Multiple Hamiltonian matrices can hence be generated
    as needed from a trajectory file (the step-by-step evolution of the position of the atomic coordinates).
    In order to generate the trajectories we run several QMD simulations using the LATTE code. 
    
    In the initial stages of our research, we chose to generate small Hamiltonian matrices of size $5\times 5$ as a proof of concept. This matrix size was selected because it allows us to efficiently test our FastEigen framework and models without the computational burden of larger matrices. The matrices were generated using the HO$^{.}$ radical system, which provides a simplified yet representative model system for our study. The results obtained from these small-scale simulations were then analyzed to assess the validity and efficacy of our approach. If our methods prove effective at this scale, they can potentially be applied to larger, more complex systems, opening up new avenues for exploration in the field. Similarly, alphatensor~\cite{Fawzi2022-qi} demonstrates proof of concept results on small scale  matrices of similar sizes.
    
    We used the LAMMPS driver to perform the MD simulation as described in \cite{latte-lammps}. Initial NVT simulations were performed with a temperature of 300K.  Trajectories at 400K, and 500K were also generated to have more variability in the matrix elements. To ensure that all the trajectories are uncorrelated, we have used different random seeds in all the simulations. We used a 0.5 fs time-step that lead to trajectories of 5 ps total simulation timescale ($10,000$ steps). 
    
    \subsection{Jacobi Algorithm}
    As briefly mentioned in the introduction, Jacobi eigenvalue method is based on Givens rotations to systematically zero out the off-diagonal elements. In this section we will extend the explanation given in the introduction. 
    Every off-diagonal element $(i,j)$ of a matrix 
    $A$ are systematically zeroed out by applying Givens rotation matrices $G^{(i,j)}$  as follows: $A'  = (G^{(i,j)})^t A G^{(i,j)}$, where, 
    $G^{(i,j)}_{ii} = G^{(i,j)}_{jj} = \cos(\theta) = c$, $G^{(i,j)}_{ij} = -G^{(i,j)}_{ij} = \sin(\theta) = s$, and $G^{(i,j)}_{kl} = \delta_{kl}$ for every other element. The latter, defines the following set of equations: 
    
    \begin{eqnarray}    
        A'_{li} =& c A_{li} - s A_{lj},\mathrm{ \, for\,} l\ne i \mathrm{\, and \,} l \ne j \\
        A'_{lj} =& s A_{li} - c A_{lj},\mathrm{\, for\,} l\ne i \mathrm{\, and \,} l \ne k \\
        A'_{ii} =& c^2 A_{ii} + s^2 A_{jj} - c s A_{ij} \\
        A'_{jj} =& s^2 A_{ii} + c^2 A_{jj} + c s A_{ij} \\
        A'_{ij} =& (c^2 - s^2)A_{ij} + cs (A_{ii} - A_{jj}) 
        \label{eqs}
    \end{eqnarray}
    
    Since, $A'_{ij} = 0$, then $(c^2 - s^2)A_{ij} + c s(A_{ii} - A_{jj}) = 0$, from where we can solve for $\tan{(\theta)} = 2A_{ij}/(A_{jj} - A_{ii})$ to then compute $c$ and $s$ as: $c = (1+ \tan(\theta)^2)^{-1/2}$ and $s = c \tan(\theta)$.
    In conclusion, every time an element needs to be zeroed out, $c$ and $s$ values are computed and matrix $A$ is modified following the above equations. 
    At convergence, the product of all the rotation matrices gives rise to matrix $U$, while $A$ gets transformed into $D$.
    Now, two questions comes to mind: What happens to the elements that are already zero? and; which elements needs to be zeroed out first? The first question is obvious by looking at the sets of equations in \ref{eqs}. If we apply $G^{(i,j)}$ all $A'_{li}$ gets ``contaminated" with $A_{lj}$ elements, even if $A_{li}$ is 0. This means that, even if we have a sparse matrix, every iteration will create ``spurious" non-zeros where there were not before. The latter means that more rotations will need to be applied until we can guarantee that all $|A_{k\ne l}|$ are less than a certain tolerance. The answer to the second question above leads to different variations of the Jacobi algorithm which focus on which is the most appropriate sequence in which the elements needs to be zeroed out. 
    In the ``Regular Jacobi" (RJ) algorithm  the elements that are zeroed out first are the highest in absolute value: $maxA = max_{k,l}(|A_{kl}|)$. This version of the algorithm keeps on jumping from element to element looking for the next $maxA$. We will refer to this strategy as MaxElement which will become our standard for comparison. 
    Although one can demonstrate that this version of the algorithm ensures that $maxA > maxA'$, it does not guarantee to lead to the fastest convergence. In fact, it turns out that systematically zeroing out the off-diagonal elements in a cyclical way has even a better performance since there is no need to introduce a max search at every iteration. One can then specifically design a cyclical sequence to zero out the elements depending on the matrix at hand. This leads to another variant of the algorithm called the ``Cyclic Jacobi" (CJ). However, again, nothing guarantees that the proposed sequence is the most optimal one. 
    In order to prove this, one can try variations constructed by permuting the sequence generated via the RJ algorithm. It turns out that certain permutations end up with a new sequence leading to a faster convergence. In this work, we prove that, a fairly good sequence leading to a faster convergence can be learned ``on-the-fly" from the off diagonal matrix elements.
    
    \subsection{AlphaZero-based matrix diagonalization}
    The AlphaGoZero framework was developed to play the game of Go and has been successfully adapted to play many other games including Tic-Tac-Toe, Connect Four, Gomoku and Chess.   
    \begin{table}[ht!]
        \centering
        \begin{tabular}{m{6.0cm}|m{6.0cm}}
            
            \hline
            \rowcolor{lightgray}
            \textbf{List of necessary adaptations:} & \textbf{Solutions:}\\
            
            \hline
            Must train with input matrices rather than from empty board &  Utilize generated Hamiltonian matrices to initialize game.\\
            \hline
            Must make diagonal entries and below non selectable.
            &  Make all diagonal entries unselectable throughout the game.\\
            \hline
            Must use floats instead of integer player id (2D$\rightarrow$3D).
            & Instead of player markers at the positions, use floating point values.\\
            \hline
            Rather than just updating a marker position, the status of the full matrix for each Jacobi rotation must change with each step.
            & Implement Jacobi rotate algorithm to act on the selected positions and update the board to generate new selectable positions.\\
            \hline
            Must implement an entirely different stopping point for matrix diagonalization.
            & Write criterion for stopping point by checking matrices after each step.\\
            \hline
        \end{tabular}
        \caption{Implementation of Jacobi Matrix Diagonalization}
        \label{tab:feadapt}
    \end{table}
    Table \ref{tab:feadapt} describes the necessary adaptations and proposed solutions to convert the AlphaGoZero framework, initially designed for board games, into a matrix diagonalizer. Key modifications include:
    
    \begin{itemize}
        \item\textbf{ Training with Input Matrices:} Instead of starting from an empty board, the adapted model uses Hamiltonian matrices to initiate the game.
        \item \textbf{Selection Restriction:} The model ensures only elements above the diagonal of a matrix are selectable throughout the game, as opposed to full board accessibility in the original version.
        \item \textbf{Data Type Transition:} The model replaces integer player IDs typically used in games with floating-point values, facilitating a transition from a 2D game board to a 3D matrix representation  (The three dimension corresponds to the matrix $\mat{S}$ co-ordinates $(x,y)$ and corresponding value ($S(x,y)$)).
        \item \textbf{Matrix Status Evolution: }The model incorporates the Jacobi rotation algorithm, which operates on selected matrix positions and updates the entire matrix status for every rotation. This contrasts with the original framework, which only updated a single position marker.
        \item \textbf{Stopping Criterion for Matrix Diagonalization:} A new stopping criterion is implemented to halt the diagonalization process. The criterion is based on the evaluation of matrices after each iteration, which differs from typical game-ending conditions in the AlphaGoZero framework.
        
    \end{itemize}
    
    These changes are part of a broader suite of modifications that also involve optimizing parameters for efficient training and inference processes.

    \begin{figure}[hb!]
        \centering
        \includegraphics[width=0.6\linewidth]{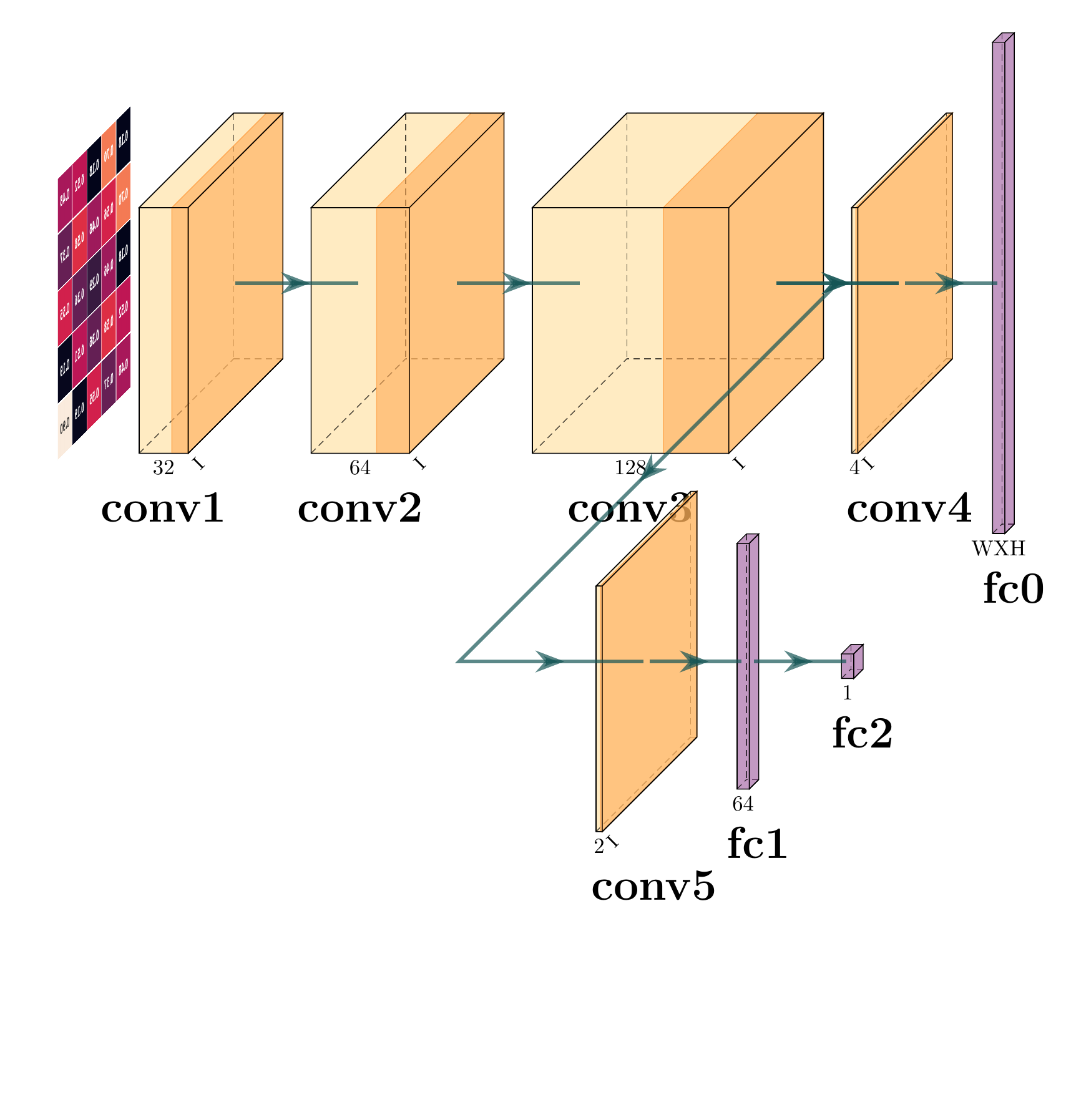}
        \caption{Overview of the policy-value Network (PVN). PVN takes a dense matrix as input and outputs a policy vector  $p_i$ and value $v_i$ representing distribution over diagonalizing steps and current player winning probability respectively. The framework includes a shared feature extraction composed of three convolutional layers (\textbf{conv1}, \textbf{conv2}, and \textbf{conv3}) used by both the policy and value network. The policy network includes an additional convolutional layer (\textbf{conv4}) and a fully connected layer (\textbf{fc0}) with output size $W\times H$, while the value network combines the shared feature extraction with \textbf{conv5}, \textbf{fc1}, and \textbf{fc2} to produce a scalar output. }
        \label{fig:dqn}
    \end{figure}
    
    The complexity requirement of traditional Max-element Jacobi rotation based matrix diagonalization approach grows cubically with the matrix size and may be infeasible for matrices beyond a certain size. To address this concern, we aim to utilize AlphaZero framework for finding the near optimal solution \cite{AlphaZero-dm}. In this work we have repurposed AlphaZero to allow reformulating matrix diagonalization as a gameplay and learn the optimal diagonalization procedure as a game-winning pattern from numerous simulations. This procedure provides an ability to scale the framework for large-sized matrices without the requirement of significantly large computation resources. 
    This is possible due to the unique combination of the MCTS framework and the policy-value network  utilized in the AlphaZero framework. 
    \begin{figure}[hb!]
        \centering
        \includegraphics[width=0.6\linewidth]{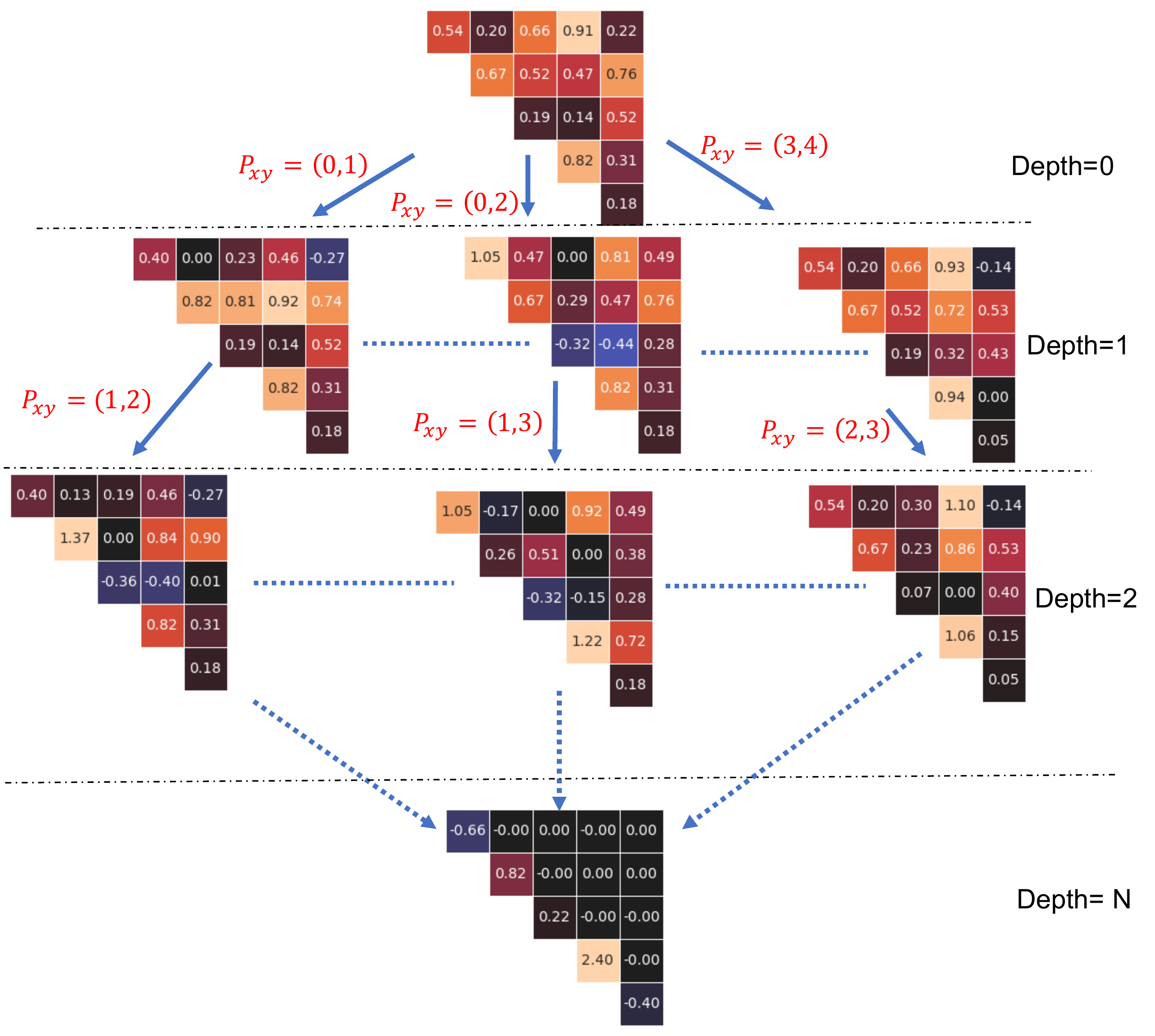}
        \caption{Overview of the MCTS graph for matrix diagonalization. The nodes in the graph represent the states of the matrix diagonalization state search where the edges $P_{xy}$ correspond to matrix rotation action. UCT formula is used to balance the exploration and exploitation. The search terminates when all off-diagonal elements are zeroed out. Depth 0 correspond to the initial matrix, depth N corresponds to the final solution, and other depths represent intermediate state transitions.}
        \label{fig:mcts}
    \end{figure}
    The policy network (upper output head) shown in Figure~\ref{fig:dqn} provides the likelihood probability for performing the Jacobi rotation corresponding to each off-diagonal elements from the upper triangular portion of the matrix (i.e., dictating what the next efficient move for performing Jacobi rotation is). In contrast, the value network (lower output head) assigns a reward  or penalty corresponding to whether the matrix was finally diagonalized or not, within some given convergence tolerance. For all the calculations throughout this work we have set the tolerance to 1e-5.
    
    The policy-value network is also combined with the MCTS to provide a lookahead search, which, when combined with the policy network, can narrow down the search to high-likelihood moves and, with the value network, can evaluate the positions in the tree through reward. An example of the tree utilized by the MCTS framework is shown in Figure~\ref{fig:mcts}.
    
    This section will discuss how we adapt the AlphaZero frameworks to solve the matrix diagonalization problem. We call it FastEigen. The beauty of the AlphaZero framework is the ability to ultimately learn the game from scratch in an unsupervised fashion and with no-domain knowledge other than the game rules. The problem of matrix diagonalization is first framed as a two-player game where the agents estimate the index for performing Jacobi rotations and compete with each other to find the fastest diagonalization. To define the 
    matrix diagonalization as a RL problem, we have following configurations:
    \begin{itemize}
        \item \textbf{States:} The current state of the  matrix being diagonalized. 
        \item \textbf{Rewards:} The reward is given based on the quality of the diagonalization. It is related to the closeness of off-diagonal elements to zero, indicating successful diagonalization.
        \item \textbf{Actions:} An action is the selection of the pivot point for the next Jacobi rotation. This is the decision that the AlphaZero model has to make at each step. In this specific case, the actions search space is constrained to off-diagonal elements from upper triangular portion of the matrix. 
        
    \end{itemize}
    
    \begin{figure}[ht!]
        \centering
        \includegraphics[width=0.9\linewidth]{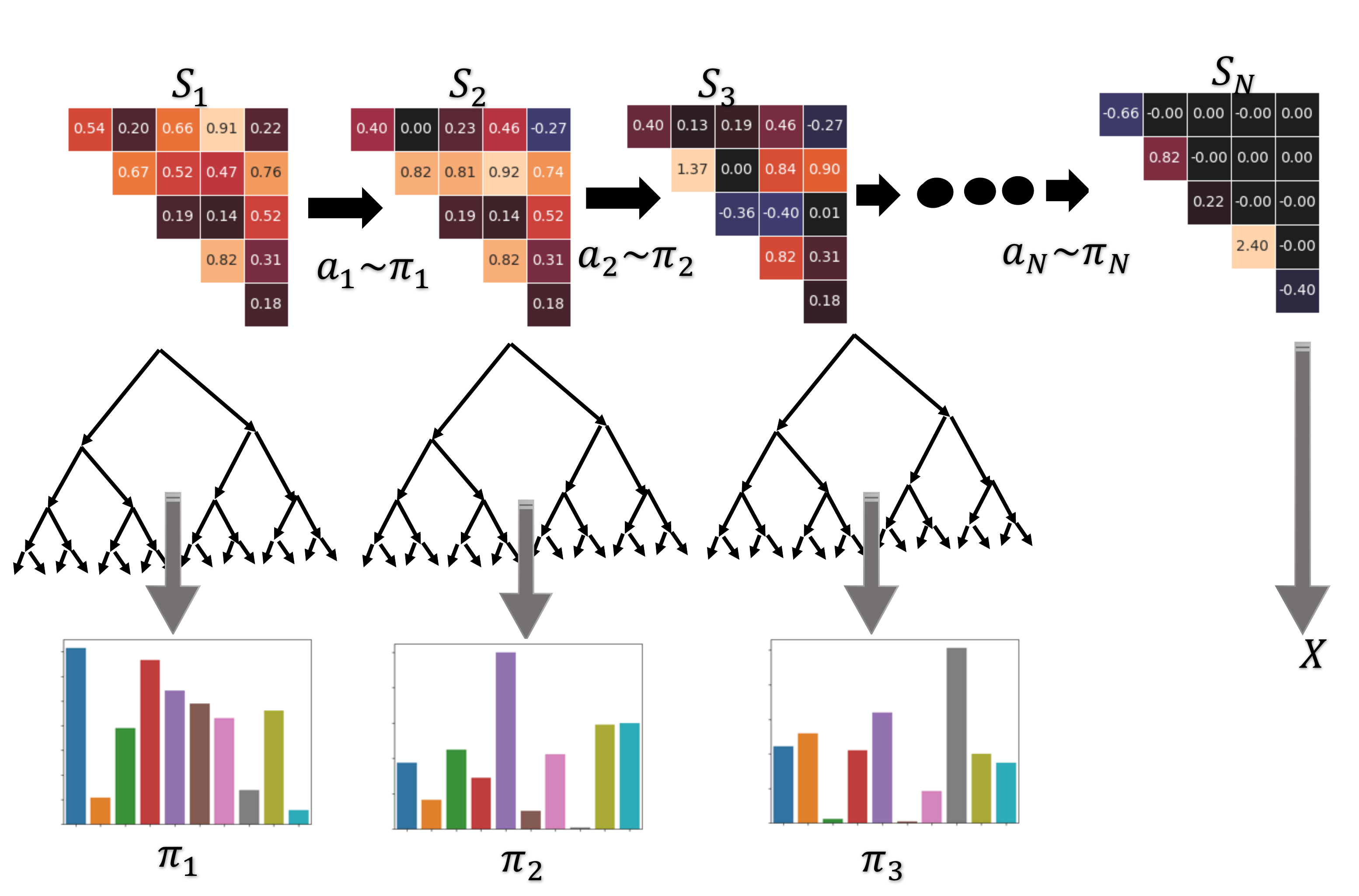}
        \caption{Overview of the self play using the MCTS framework for matrix diagonalization. During self-play, the agent self evaluates a sequence of matrix diagonalizations $S_1,S_2,...S_N$ where an MCTS parameterized by $\alpha_{\theta}$  is executed at each stage with latest policy-value network parameterized by $f_{\theta}$. The next diagonalization step is based on the policy distribution estimated by MCTS($a_i\sim \pi_i$). $X$ corresponds to the terminal diagonalization state when all the off-diagonal elements have been zeroed out.} 
        \label{fig:mcts_policy}
    \end{figure}
    The FastEigen system comprises two major building blocks (i) the policy-value network as shown in Figure~\ref{fig:dqn} \cite{Silver2013-se,Watkins_Christopher_undated-kh} and ii) MCTS as shown in Figure~\ref{fig:mcts}. The role of the policy value network is to observe the current matrix state and generate a decision for the next rotation operation to be performed to diagonalize the matrix with the least number of steps. On the other hand, MCTS contemplates the multiple possible outcomes starting from the current one. The MCTS simulation can be expressed in the form of the tree shown in Figure~\ref{fig:mcts}. MCTS then  utilizes the decision outcomes from the policy-value network to control the simulation procedure where the output of the policy value network is shown in Figure~\ref{fig:mcts_policy}. Similarly, the policy value network is trained on the simulations of MCTS as shown in Figure~\ref{fig:mcts_train}. The policy value framework is a deep learning framework that takes the matrix state and outputs the probability over the actions and the winning chance of the current state. The parametric representation is given as follows:
    \begin{equation}
        (\vect{p},\mat{v}) = f_{\theta}(s)
    \end{equation}
    where $ f_{\theta}(s)$  represents the neural network with parameter $\theta$. For a given matrix of a size $n\times n$, the state matrix $\mat{s}$ is the $n\times n$ array representing the  matrix. 
    
    The Jacobi rotation method is an algorithm specifically designed for symmetric matrices. 
    Because of this symmetry, the action space for AlphaZero is constrained to only the upper triangle elements of the matrix. 
    Here $\vect{p} = (p_1,p_2,...,p_N)$ where $N = n(n-1)/2$ is the policy distribution vector whose $i^{th}$ element $p_i = P_r(a|s)$ correspond to the prior probability of performing Jacobi rotation on the $i^{th}$ element. Here, $a$ corresponds to the most probable action based on $\vect{p}$ distribution. Similarly, $v\in [-1,1]$  represents the winning value for the current player who will perform the next Jacobi rotation. Larger $v$ corresponds to a higher chance of being able to diagonalize.
    
    A deep neural network with convolutional and fully connected layers was utilized for the policy value network as shown in Figure~\ref{fig:dqn} with two different output heads for the policy and the value outputs. This framework comprises the following module:
    \begin{figure}[ht!]
        \centering
        \includegraphics[width=0.9\linewidth]{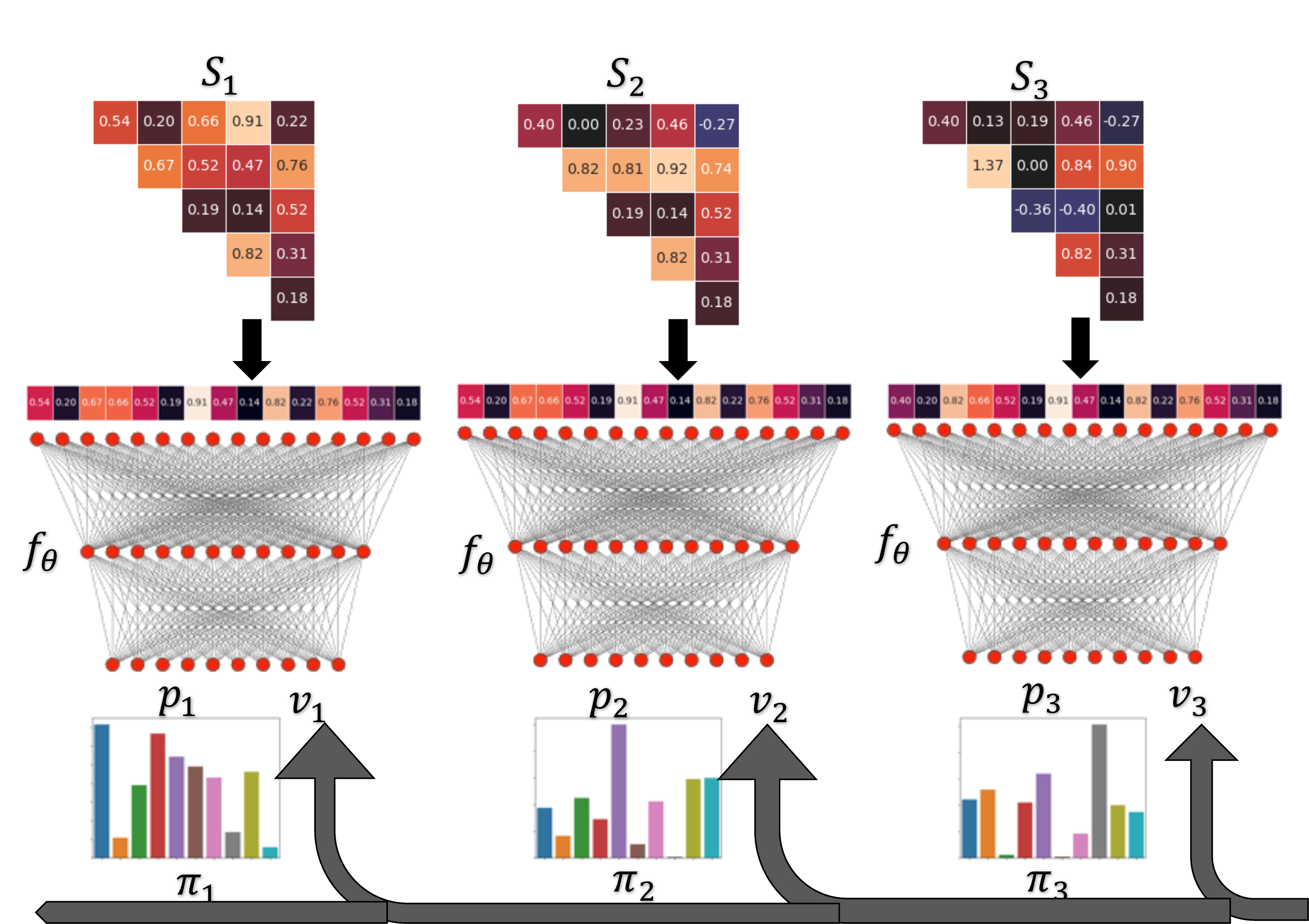}
        \caption{Overview of training of the FastEigen framework. Here $S_i, f_\theta, p_i$ and $v_i$ correspond to State, neural network parameter, policy distribution, and predicted current winner winning probability value for $i^{th}$ game respectively. The bar graph illustrates $p_i$, with distinct color coding corresponding to each diagonalization pivot's associated probability, from which the action with the highest probability is selected as the next pivot for the Jacobi rotation. During the training, the NN parameters $f_{\theta}$ is updated so as to maximize $p_i$ and MCTS search probabilites $\pi_i$ similarity and minimize error between $v_i$ and game winner X. }
        \label{fig:mcts_train}
    \end{figure}
    The policy value framework comprises of  \textbf{Feature extraction module:} This module takes the input matrix as the state space and extracts the features. These features are then fed separately to the policy and value head,
    \textbf{Policy head:} This module generates the prior probability distribution vector $\vect{p}$.
    and \textbf{Value head:} This module generates the winning value $v$.
    
    The MCTS framework is collectively parameterized by  $\alpha_{\theta}$ and it is instructed by the policy-value neural network collectively parameterized by $f_{\theta}$. The policy value network receives the current state matrix $\mat{s}$ and outputs a policy probability $\vect{\pi}$ which correspond to the diagonalization pivot's associated probability.  
    The policy distribution vector, depicted as $\pi = (\pi_1, \pi_2, ..., \pi_{N})$ in Figure~\ref{fig:mcts_policy}, represents the probabilities associated with executing matrix diagonalization for each respective off-diagonal pivot. The network selects the pivot associated with the highest probability value for the next operation.
    
    In the FastEigen framework, the output of the MCTS, represented by $\alpha_{\theta}$, is given by $\pi = \alpha_{\theta}(s)$. Here, $\pi$ is the policy distribution vector, representing the probabilities of selecting each off-diagonal pivot for matrix diagonalization. These probabilities are influenced by a temperature parameter $\tau$, which controls the balance between exploration and exploitation. More specifically, each action's probability $\pi_a$ is proportional to $N(s, a)^{1/\tau}$, where $N(s, a)$ is the number of times action $a$ has been taken in state $s$. Higher $\tau$ values make the probabilities more uniform, encouraging exploration of different pivots, while lower $\tau$ values make the probabilities more skewed towards the most frequently selected pivots, promoting exploitation.
    
    In the MCTS framework, the search tree  comprises a graph $(s,a)$ where the edge stores the prior proability $P(s,a)$, node visit count  $N(s,a)$ and action-value $Q(s,a)$. The simulation starts with the root node and iteratively selcts the moves that maximize the upper confidence bound $\zeta =  Q(s,a) + U(s,a)$ where $U(s,a) \propto P(s,a)/(1+N(s,a))$ imply $U(s,a) = C_{puct} P(s,a)/(1+N(s,a)) $ and $\zeta =  Q(s,a)+C_{puct} P(s,a)/(1+N(s,a))$ where $C_{puct}$ is a proportional constant that trades between exploitation term (first) and exploration term (second), until the leaf node $s'$ is encountered. The leaf node is then finally evaluated by the  policy-value network $f_{\theta}$ to generate prior probabilities and evaluation as $(p(s'),v(s')) = f_{\theta}(s')$. During the game play, each edge $(s,a)$ traversed is updated to increment the visit count $N(s,a)$ and the action-value is updated to the mean as $Q(s,a) = \frac{1}{N(s,a)} \sum_{s'|s,a \rightarrow s'} v(s')$ where $s,a \rightarrow s'$ corresponds to reaching to state $s'$ from state $s$ after taking action $a$.
    
    To train the policy value framework, first, MCTS is deployed to play each move as per the self-play nature of the RL employed as shown in Figure~\ref{fig:mcts_train}. The neural network $f(\theta)$ is first initialized randomly so that the initial weights are $\theta_0$. Then for each subsequent iteration $i \geq 1$, games of self-play are generated, and for each time step $t$, the MCTS search policy is estimated as $\pi_t = (\alpha_{\theta})_{i-1}(s_t)$ using a neural network with parameters corresponding to the previous iteration, i.e., $f_{{\theta}_{i-1}}$ where the moves are then played by sampling from search probabilities $\pi_t$. At each time step $t$, the data is stored as $(s_t,\pi_t,z_t)$ where $z_t = \pm r_T$ is the winner of the game and new network parameters $\theta_i$ are trained from data $(s,\pi,z)$ sampled uniformly among all time steps of the last iteration of self-play. $r_T$ represents the final reward of the game, given at the terminal state $T$. Here, the final reward $r_T$ could be +1 for a win, -1 for a loss, and 0 for a draw. Therefore, $z_t$ being equal to $\pm r_T$ would mean that $z_t$ takes the value of the final game reward, with its sign indicating whether the outcome was a win or a loss. The policy value network $f_{\theta}$ is adjusted to minimize the error between the predicted value $v$ and the self-play $q$ index $z$ and to maximize the similarity of neural network move probabilities $p$ to the search probabilities $\pi$ where the parameters $\theta$ are adjusted with gradient decent on loss function $l$ which is given as
    \begin{equation}
        l = (z-v)^2 - \pi^T \log(p) + c ||\theta||^2
    \end{equation}
    The first term in the loss function is the Mean Square Error (MSE) error, the second term is the cross entry, and the third term is the $l2$-regularizer term to prevent overfitting.
    
    A comprehensive overview of the algorithm is presented in the pseudocode of Algorithm \ref{alg1}. 
    
    \begin{algorithm}[tb]
        \caption{AlphaFastEigen framework}
        \label{alg1}
        \begin{algorithmic}
            \STATE {\bfseries Input:} data $X$, $\theta$
            \STATE Initialize DNN $f_{\theta}$
            \REPEAT
            \STATE Play Game
            \STATE Update $\theta$ 
            \WHILE{$\neg$ Win $\lor$ $\neg$ Loss}\STATE \% Until win/loss continue play
            
            \STATE From current state s, perform MCTS
            \STATE Estimate move probabilities $\pi$ by MCTS
            \STATE Record $(s_{\pi})$ as an example
            \STATE Randomly draw the next move from $\pi$ 
            \ENDWHILE
            \STATE Update  $\theta$
            \STATE Let z be previous game outcome $(+1 \, \textup{or} \, -1)$
            \STATE Sample from last game's examples $(s, \pi, z)$ 
            \STATE Train DNN $f_{\theta}$ on sample to get new $\theta$
            \STATE \textbf{repeat until} no change in  $\theta$
            \UNTIL{no change in $\theta$}
        \end{algorithmic}
    \end{algorithm}
    
    \section{Experiments and Results}
    
    Three trajectories for HO$^{\cdot}$ were created and processed into 1000 $5\times 5$ matrices sequences.  The temperatures for each sequence were 300K, 400K and 500K.  Each of these matrices sequences was then split randomly into a 750 matrices set for training and a 250 matrices set for testing.  Training was conducted with the exploration-exploitation parameter $C_{puct}$ set to 4, MCTS playouts parameter $n_{playouts}$ set to 560 and the action policy network was trained for 50 epochs. $C_{puct}$ is a hyperparameter that balances exploration and exploitation in the MCTS algorithm where value of 0 corresponds to pure exploitation and larger values favor exploration over exploitation.$n_{playouts}$ corresponds to the number of simulations or ``playouts" run from each position during the MCTS process.The optimal value for $n_{playouts}$ represents a trade-off: higher values can potentially improve output performance but at the cost of increased training time. The training time was 5.8 hours on a single NVIDIA A100 GPU with an AMD EPYC 7713 64-Core cpu for $5\times 5$ sized matrices.  
    \begin{figure*}[!ht]
        \centering
        \includegraphics[width=1\linewidth]{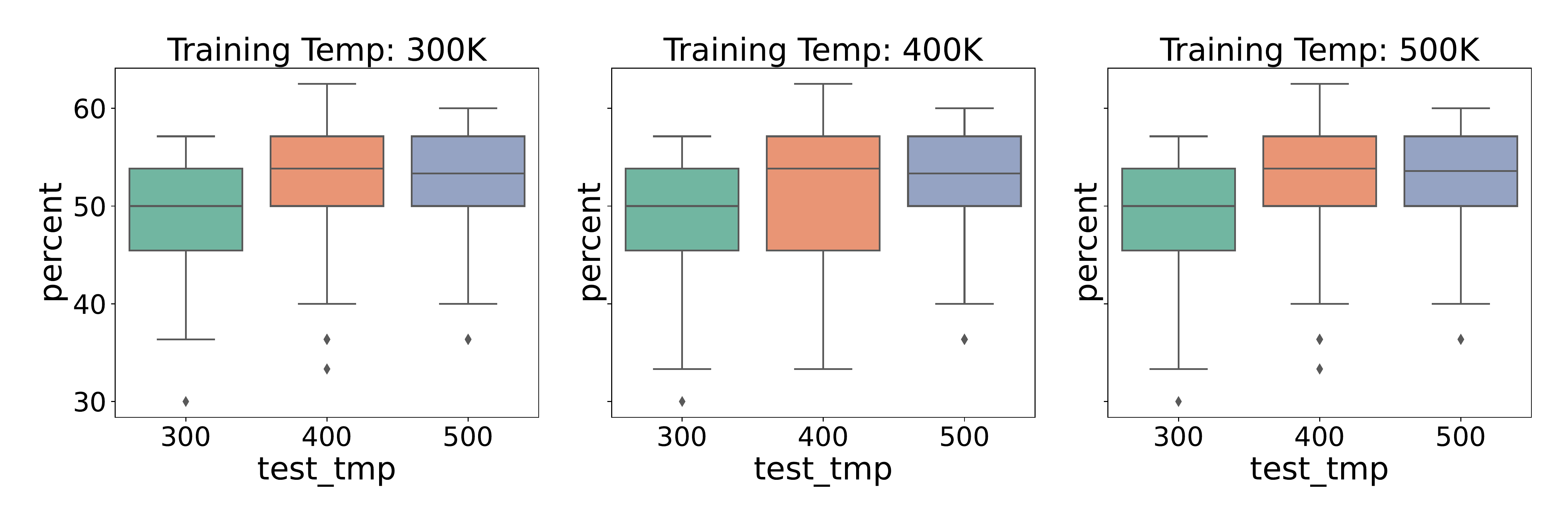}
        \caption{Overview of net savings with AlphaFastEigen compared with MaxElem based matrix diagonalization}
        \label{fig:net_savings}
    \end{figure*}
    
    Our inference was carried out on nodes outfitted with 4 NVIDIA A100 GPUs and powered by AMD EPYC 7713 64-Core CPUs. We employed 12 independent GPU processes (MPS servers) for each node in a embarrassingly parallel  computing approach. We capped the inference time at 5 minutes per matrix to create all solution paths. This translated to a maximum of roughly 10 hours for the complete inference phase per set of 1000 matrices.
    \begin{figure*}[!ht]
        \centering
        \includegraphics[width=0.9\linewidth]{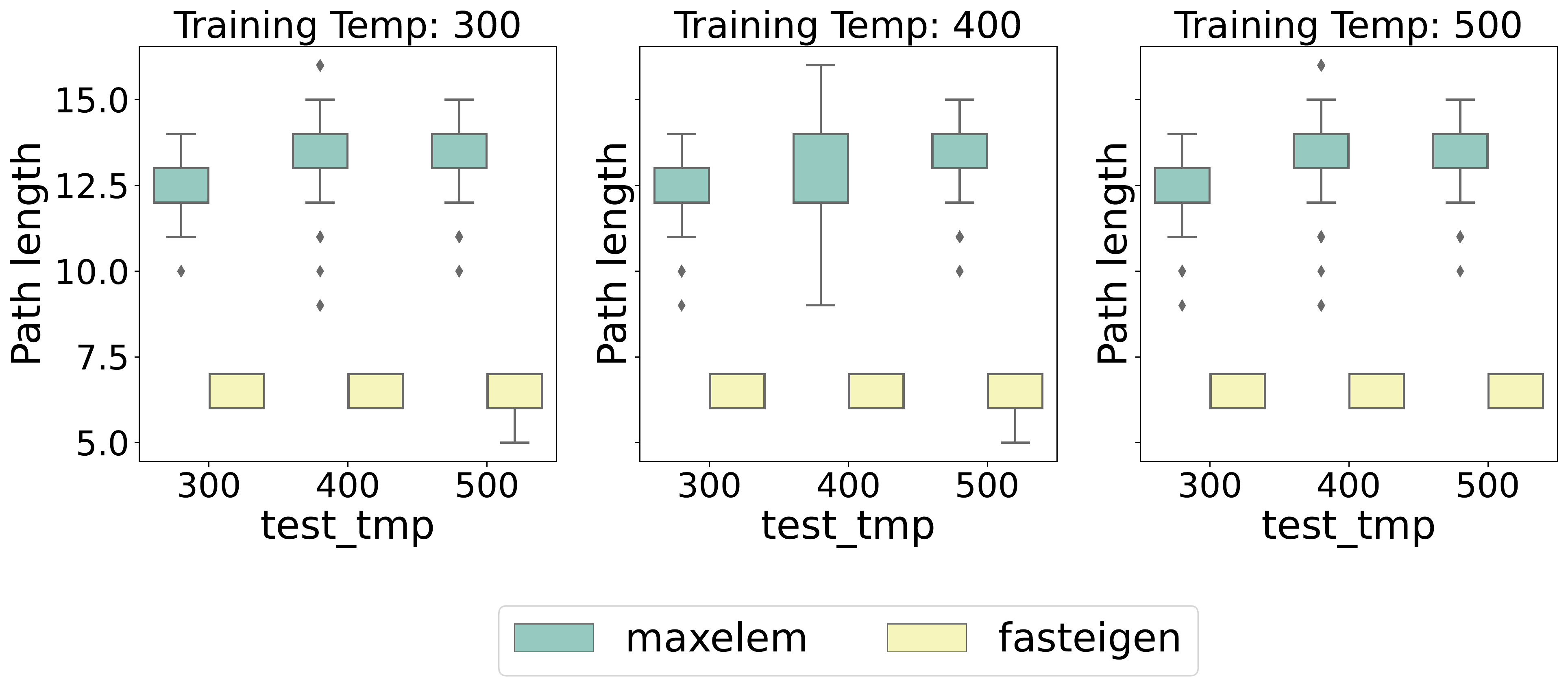}
        \caption{Overview of overall path length based on AlphaFastEigen compared with MaxElem based matrix diagonalization}
        \vspace{1cm}
        \label{fig:fasteigenvsmaxelem}
    \end{figure*}
    For optimal test performance, the inference experiments involved considerable tuning  of the $C_{puct}$ and $n_{playouts}$ parameters to produce solutions with the FastEigen system.   Unlike finding the patterns for playing the board Games (which AlphaZero was initially designed for), finding the optimal path for the matrix diagonalization is a challenging task where different matrices might have different underlying diagonalization patterns. It is challenging for AlphaZero to find a one-shot solution for diagonalization steps that work for all matrices. Due to this challenge, our model requires significant tuning associated with choosing the right set of hyper-parameters during the inference. Due to this fact, the model required significant exploration steps along with exploitation, even during inference, which is uncommon for many RL frameworks. Moreover, due to the nature of the matrix diagonalization, even for a smaller matrix, many paths could lead to the same solution with FastEigen. Although most matrices can produce efficient solutions in seconds, a small minority can continue producing solutions for several hours due to the exploration of inefficient paths in MCTS. Setting the parameters of $c_{punct}$ and $n_{playouts}$ to 4 and 13,000 respectively, produced the results presented here, and inferencing times were limited to a maximum of 5 minutes for each matrix.
    
    \begin{figure*}[!ht]
        \centering
        \includegraphics[width=0.99\linewidth]{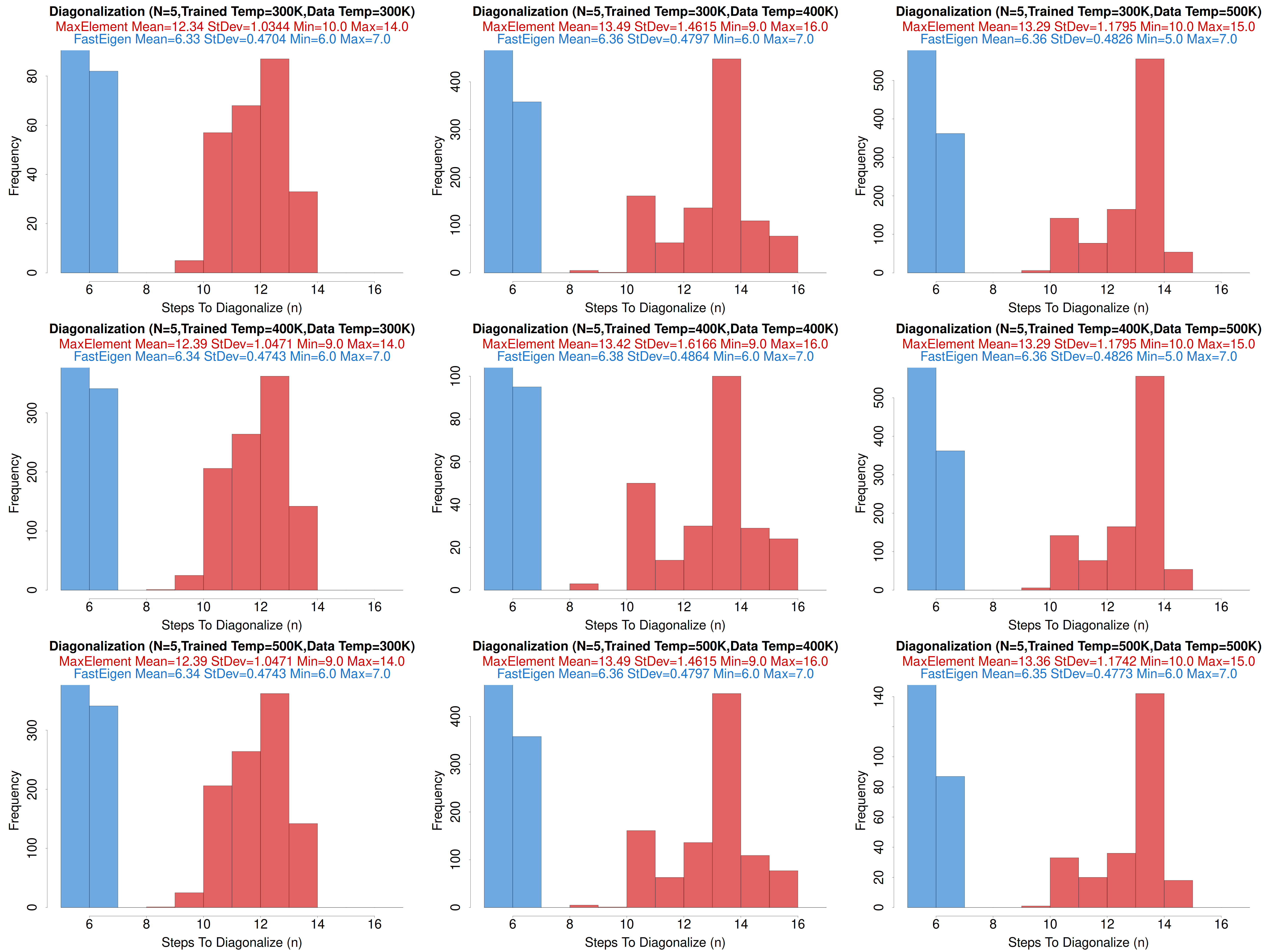}
        \caption{Overview of path length histogram based on FastEigen compared with MaxElem based matrix diagonalization}
        \label{fig:fasteigenvsmaxelem_hist}
    \end{figure*}
    In each case, the FastEigen method consistently outperformed, delivering solutions within six to seven steps for all 250 test matrices achieving large savings ($>50\%$) in path length as shown in Figures~\ref{fig:net_savings}. On the other hand, the MaxElement method required a broader range of steps, varying from 10 to 14 at 300K, 11 to 16 at 400K, and 10 to 15 at 500K as shown in Figures~\ref{fig:fasteigenvsmaxelem}. While these results underscore the noticeable superiority of the FastEigen method over the MaxElement method, they do not provide an insight into the distribution of solution advantages. The subsequent histogram plots will provide a clearer depiction of how these solution advantages are distributed shown in Figures~\ref{fig:fasteigenvsmaxelem_hist}.
    
    \begin{figure*}[!ht]
        \centering
        \includegraphics[width=0.99\linewidth]{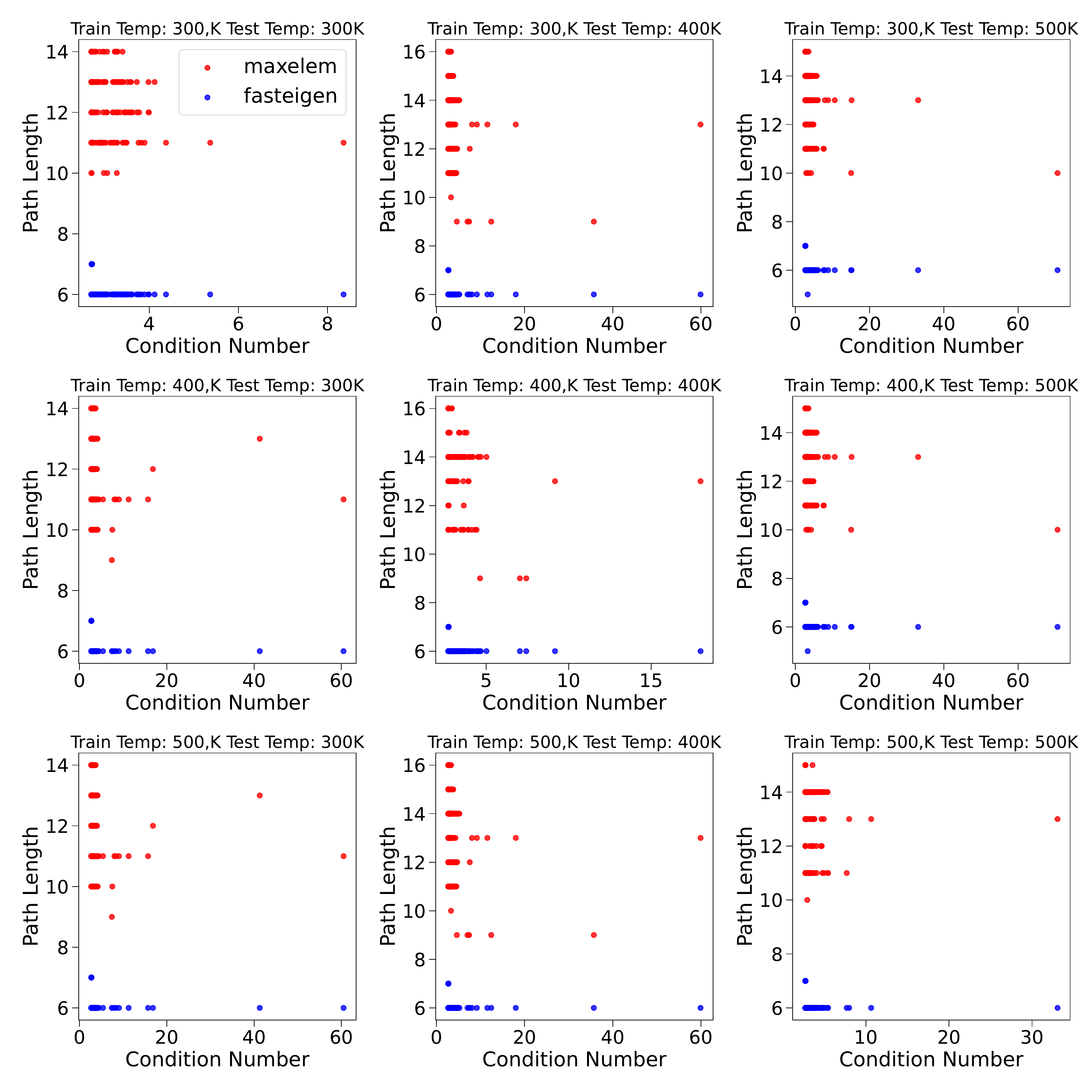}
        \caption{Correlation between condition number and path lengths corresponding to FastEigen compared with MaxElement}
        \label{fig:fasteigenvsmaxelem_cond}
    \end{figure*}
    
    The results show that the FastEigen method allows for a significant efficiency improvement in matrix diagonalization. These improvements were observed in matrices where the MaxElement method had an average of 12.36 steps to solve. The FastEigen method demonstrated a substantial reduction in the average number of steps, ranging from 6.05 to 7.01. It's also noteworthy that the FastEigen method provided a substantial advantage over the MaxElement method across all tested matrices. Furthermore, we examined the adaptability of the FastEigen method across different temperatures by running inference on the complete sets of 1000 matrices that were not trained at the same temperature. The results from these mixed temperature solutions for each test dataset are presented in the plots of Figure~\ref{fig:fasteigenvsmaxelem_hist}.

    Utilizing the FastEigen method instead of the MaxElement method yielded substantial average savings, ranging between 48.5\% and 52.3\% (Figure~\ref{fig:net_savings}). The limited variability in these percentages suggests that the FastEigen method is highly adaptable across different temperatures. This implies that there is no need to recalibrate the FastEigen method for each specific temperature to achieve substantial efficiency gains. Additionally, we have illustrated the distribution of the percentage savings achieved by the FastEigen method over the MaxElement method for each cross-temperature variation in the Figures~\ref{fig:fasteigenvsmaxelem_hist}.

    Based on the resulting plots, there are no large discrepancies when training on one temperature and inferencing on another indicating that training is fairly robust and need not be changed to inference on other temperatures.

    In the comparison of the MaxElement and FastEigen methods for matrix diagonalization, we observed distinct behaviors between the two. Particularly, the number of steps required for MaxElement exhibited significant variation across different matrices, even when those matrices shared the same condition number as shown in Figure~\ref{fig:fasteigenvsmaxelem_cond}. This inconsistency could potentially introduce unpredictability when implementing the MaxElement method in practical scenarios, as the computation cost might vary widely across different cases. On the other hand, Fasteigen demonstrated a much more consistent performance. The number of steps required by FastEigen remained relatively stable across matrices with different temperatures and condition numbers. This consistency suggests that FastEigen provides a more reliable estimate of the computation cost involved in matrix diagonalization, which can be a valuable attribute in many practical applications. Furthermore, the consistent behavior of FastEigen implies that it may be less sensitive to variations in the properties of the input matrices, offering a more robust solution for matrix diagonalization across a wide range of scenarios.

    \section{Conclusions}
    
    By formulating a matrix diagonalization as a board game we have demonstrated how reinforcement learning using the AlphaZero AI framework can be used to learn the fastest path to solution.  A significant acceleration can be demonstrated for the FastEigen method. 
    The accelerated performance of the FastEigen method over the MaxElement method was consistently highlighted, with substantial savings observed in the number of steps needed to solve matrices sampled from QMD simulation trajectories across various temperatures. This robust performance of the FastEigen method, despite training on one temperature and inferencing on another, signifies its adaptability and resilience, making it a promising tool for wide-ranging applications.
    These findings emphasize the potential of leveraging advanced reinforcement learning techniques in matrix diagonalization, with significant improvements over established methods. 
    Our findings highlight the opportunity to use machine learning as a promising tool to improve the performance of numerical linear algebra.
    Looking ahead, our research endeavors will concentrate on extending this approach to handle larger Hamiltonian matrices. The promising results from this study suggest a good foundation for further optimization and scaling up, with the ultimate goal of advancing the capabilities of matrix diagonalization methods in computational physics and other similar fields.
    
      \section*{Acknowledgments}
  This manuscript has been approved for unlimited release and has been assigned LA-UR-23-21573. 
    This work was supported by the Laboratory Directed Research and Development program of Los Alamos National Laboratory under project number 20220428ER.   This research used resources provided by the Los Alamos National Laboratory Institutional Computing Program. Los Alamos National Laboratory is operated by Triad National Security, LLC, for the National Nuclear Security Administration of U.S. Department of Energy (Contract No. 89233218CNA000001). Additionally, we thank the CCS-7 group, Darwin cluster, and Institutional Computing (IC) at Los Alamos National Laboratory for computational resources. Darwin is funded by the Computational Systems and Software Environments (CSSE) subprogram of LANL’s ASC program (NNSA/DOE).

    \bibliography{biblio}
    \bibliographystyle{unsrt}
    
\end{document}